\documentclass[reqno]{amsart}

\usepackage{amsthm,amsmath,amssymb,amsfonts,bm}
\usepackage{stmaryrd}
\usepackage{bbold}
\usepackage{times}
\usepackage[all]{xy}
\usepackage{mathrsfs}
\usepackage{graphicx}
\usepackage{tikz}
\usetikzlibrary{matrix,calc,arrows}

\theoremstyle{plain}

\newtheorem*{ST}{Stratification Theorem}

\theoremstyle{definition}

\theoremstyle{remark}

\newcommand{\Ker}{\operatorname{Ker}}

\newcommand{\Aut}{\operatorname{Aut}}

\newcommand{\Spec}{\operatorname{Spec}}
\newcommand{\Fract}{\operatorname{Fract}}
\newcommand{\Supp}{\operatorname{Supp}}

\newcommand{\tto}{\longrightarrow}
\newcommand{\iso}{\stackrel{\sim}{\tto}}

\renewcommand{\k}{\mathbb{k}}
\renewcommand{\phi}{\varphi}

\newcommand{\GSpec}{G\text{-}\!\Spec}

\newcommand{\cC}{\mathcal{C}}

\newcommand{\fp}{\mathfrak{p}}

\newcommand{\fa}{\mathfrak{a}}
\newcommand{\fb}{\mathfrak{b}}
\newcommand{\fq}{\mathfrak{q}}

\newcommand{\e}{\varepsilon}

\newcommand{\ZZ}{\mathbb{Z}}

\newcommand{\bx}{\mathbf{x}}

\newcommand{\byG}{\!\!:\!\! G}

\newdir{ >}{{}*!/-5pt/\dir{>}}

\begin{document}




\title[Stratification of Noncommutative Prime Spectra]%
{On the Stratification of Noncommutative Prime Spectra}

\author{Martin Lorenz}

\address{Department of Mathematics, Temple University,
    Philadelphia, PA 19122}

\email{lorenz@temple.edu}

\urladdr{http://www.math.temple.edu/$\stackrel{\sim}{\phantom{.}}$lorenz}


\thanks{Research of the author supported in part by NSA Grant H98230-12-1-0221}



\maketitle

\begin{abstract}
We study rational actions of an algebraic torus $G$ by automorphisms on an
associative algebra $R$. The $G$-action on $R$ induces a stratification of the
prime spectrum $\Spec R$ which was introduced by Goodearl and Letzter.
For a noetherian algebra $R$, Goodearl and Letzter showed that
the strata of $\Spec R$ are isomorphic to the spectra of certain 
commutative Laurent polynomial algebras. The purpose of this note is to give
a new proof of this result which works for arbitrary algebras $R$.
\end{abstract}

\maketitle



\section*{Introduction}

Let $G$ be an affine algebraic group and let $R$ be an associative
algebra on which $G$ acts rationally by algebra automorphisms. 
The induced $G$-action on the set $\Spec R$ of all prime
ideals of $R$ leads to a stratification of $\Spec R$ which was pioneered by
Goodearl and Letzter \cite{kGeL00}. For the special case of an 
algebraic torus $G$ and a noetherian algebra $R$, Goodearl and Letzter
have given a description of the strata of $\Spec R$ in terms of the
spectra of certain (commutative) Laurent polynomial algebras. 
Later, a different description of the strata was given in \cite{mL09}, for any algebra $R$
and any connected affine algebraic group $G$.
The purpose of this short note is to consolidate the description of \cite{mL09}, for the case of an
algebraic torus $G$,
with the earlier one due to Goodearl and Letzter, still working with a general algebra $R$.

Throughout, $\k$ will be an algebraically closed base field of
arbitrary characteristic. The largest $G$-stable ideal of the algebra $R$ that is contained in 
a given ideal $I$ of $R$, called the \emph{$G$-core} of $I$, will be 
denoted by $I\byG$; so $I\byG = \bigcap_{g \in G} g.I$.
It is easy to see \cite[Proposition 8(b)]{mL08} that, for any rational action of
an affine algebraic group $G$ on $R$,  the collection of all 
$G$-cores of prime ideals of $R$
coincides with the set of all $G$-prime ideals of $R$; the latter set
will be denoted by $\GSpec R$. If the algebraic group $G$ is connected, 
then $\GSpec R$ is simply the set of all $G$-stable prime ideals of $R$
\cite[Proposition 19(a)]{mL08}.
The \emph{Goodearl-Letzter stratification} of $\Spec R$ is the partition
\begin{equation*} 
\Spec R = \bigsqcup_{I \in \GSpec R} \Spec_IR 
\qquad \text{with} \qquad
\Spec_IR = \{ P \in \Spec R \mid P\byG = I \}\ .
\end{equation*}
For an
algebraic torus $G$, we will give a description of each stratum $\Spec_IR$ in terms of the spectrum of
a suitable affine \emph{commutative} algebra $Z_I$. The construction of the algebra $Z_I$ 
and the precise statement of the main result will be given in 
Section~\ref{S:statement} while the proof will occupy Section~\ref{S:proof}.


\section{The Stratification Theorem} \label{S:statement}

Let $G$ be a connected affine algebraic group over $\k$
and let $R$ be an associative
$\k$-algebra with a rational $G$-action by $\k$-algebra
automorphism.
For a given $I \in \GSpec R$, let $\cC(R/I)$ denote the extended centroid of 
the algebra $R/I$; this is a $\k$-field, called the \emph{heart of $I$}, 
on which $G$ acts via
its action on $R/I$ \cite[2.3]{mL08}. We put
\[
Z_I = \{ c \in \cC(R/I) \mid \text{the orbit $G.c$ spans a 
finite-dimensional $\k$-subspace of $\cC(R/I)$}\}\ .
\]
Clearly, $Z_I$ is a $G$-stable $\k$-subalgebra of $\cC(R/I)$ and $Z_I$ contains
the subfield of $G$-invariants, $\cC(R/I)^G$. Moreover, the $G$-action on
$Z_I$ is rational by \cite[Lemma 18(b)]{mL08}.

We now focus on the case of an algebraic torus $G$.
As usual, $X(G)$ will denote the lattice of rational characters of $G$.
The following theorem, under the additional assumption that the algebra $R$ 
is noetherian,  is originally due to Goodearl and Letzter
\cite{kGeL00}; see also \cite{kG00}, \cite{kGjtS00} and \cite[II.2.13]{kBkG02}. 

\begin{ST}
Let $G$ be 
an algebraic torus over $\k$ that acts 
rationally by algebra automorphisms on the $\k$-algebra $R$,
and let $I \in \GSpec R$. Then:
\begin{enumerate}
\item
There is an isomorphism
\[
Z_I \cong \cC(R/I)^G\Gamma_I\ ,
\] 
the group algebra over 
the field $\cC(R/I)^G$ of the
sublattice $\Gamma_I = X(G/\Ker_G(Z_I)) \subseteq X(G)$. 
\item
There is a $G$-equivariant order
isomorphism
\[
\gamma \colon \Spec_IR \stackrel{\sim}{\tto} \Spec(Z_I)\ .
\]
\end{enumerate}
\end{ST}


\section{Proof of the Stratification Theorem} \label{S:proof}

\subsection{Reductions and preliminaries}

We begin with some remarks that hold for an arbitrary
connected affine algebraic group $G$ over $\k$.
In order to describe the $G$-stratum $\Spec_IR$ for a given
$I \in \GSpec R$, we may replace $R$ by $R/I$ and thus assume that $I=0$. 
In particular, $R$ is a prime ring.
We will write $\cC = \cC(R)$ and $Z = Z_0$ for brevity; so $\cC$ is a commutative 
$\k$-field on which $G$ acts by automorphisms and
\[
Z =  \{ c \in \cC \mid \text{the orbit $G.c$ spans a 
finite-dimensional $\k$-subspace of $\cC$}\}\ .
\]
Our goal is to give a description of $Z$ and to
establish a suitable order isomorphism 
\[
\gamma \colon \Spec_0R 
= \{ P \in \Spec R \mid P\byG = 0 \} \stackrel{\sim}{\tto} \Spec(Z)\ .
\]
We will need  to consider various $G$-actions; they will usually be indicated by
a simple dot as in the foregoing. When more precision is
necessary, the $G$-action on $\cC$
will be denoted by $\rho$. The group $G$ also acts on 
$\k[G]$, the algebra of regular functions of $G$, via 
the right and left regular representations
$\rho_r, \rho_\ell \colon G \to \Aut_{\k\mbox{-}\text{alg}}(\k[G])$;
they are defined by $\left(\rho_r(x)f\right)(y) = f(yx)$ and
$\left(\rho_\ell(x)f\right)(y) = f(x^{-1}y)$ for $x,y \in G$ and $f \in \k[G]$.

\subsection{The algebra $Z$}

Consider the Hopf $\cC$-algebra 
\[
S = \cC \otimes_\k \k[G]\ ;
\]
this is an algebra of $\cC$-valued functions on $G$ via 
$(\sum_i c_i \otimes f_i)(g) = \sum_i c_if_i(g)$. 
The group $G$ acts on the ring $S$ via $\rho \otimes \rho_r$.
If $g \in G$ and $s = \sum_i c_i \otimes f_i \in S$ then $g.s \in S$ 
is the function $G \to \cC$ that is given by 
\[
(g.s)(x) = \sum_i (g.c_i) f_i(xg) = g.\sum_i c_i f_i(xg) \quad (x \in G)\ .
\]
Let $S^G$ denote the subring of $G$-invariants in $S$. Thus, 
$s \in S^G$ if and only if $g^{-1}.s(1) = s(g)$ for all $g\in G$.

We claim that
\begin{equation} \label{E:inZ}
S^G \subseteq Z \otimes_\k \k[G]\ .
\end{equation}
To see this, let $s = \sum_1^r c_i \otimes f_i \in S^G$ with $\{ f_i \}_1^r \subseteq \k[G]$
chosen $\k$-linearly independent. 
Choose $\{ x_j\}_i^r \subseteq G$ such that the matrix $A_g = (f_i(x_jg))_{i,j}$
is invertible and put $B = (f_i(x_j))_{i,j}$. Then the equations
$\sum_i c_if_i(x) = \sum_i (g.c_i) f_i(xg)$ for all $x,g \in G$, with $x=x_j$
for $j=1,\dots,r$, can be written as the matrix equation 
$(c_1,\dots,c_r)B = (g.c_1,\dots,g.c_r)A_g$ or else 
$(g.c_1,\dots,g.c_r) = (c_1,\dots,c_r)BA_g^{-1}$. This shows that 
$G.c_i \subseteq \sum_{i=1}^r \k c_i$ and hence $c_i \in Z$, proving \eqref{E:inZ}.

Next, we show that there is an isomorphism
\begin{equation} \label{E:Ziso1}
\e \colon S^G \iso Z 
\end{equation}
such that
\begin{equation} \label{E:intertwining}
\e \circ (1_Z \otimes \rho_l(g)) = \rho(g) \circ \e 
\end{equation}
holds for all $g \in G$. Indeed, by \cite[Lemma 18(b)]{mL08}, the $G$-action on $Z$ is rational: it arises
from a map of $\k$-algebras $\Delta_Z: Z \to Z \otimes _\k \k[G]$,
$c \mapsto \sum c_0 \otimes c_1$, via $g.c = \sum c_0c_1(g)$.
By \cite[equations (17) and (18)]{mL09}, the $\k[G]$-linear extension of $\Delta_Z$, 
which will also be denoted by
$\Delta_Z$, is an isomorphism of $\k[G]$-algebras 
\[
\Delta_Z \colon Z  \otimes _\k \k[G] \iso Z \otimes _\k \k[G]
\]
that satisfies the ``intertwining formula'' 
$\Delta_Z\circ (\rho\otimes \rho_r)(g) = (1_Z \otimes \rho_r)(g)\circ\Delta_Z$
for each $g\in G$. It follows that $\Delta_Z$ yields an isomorphism of
$S^G = (Z \otimes_\k \k[G])^G$ with the sub algebra of 
$(1_Z \otimes \rho_r)(G)$-invariants in $Z\otimes_\k \k[G]$. Since the latter
algebra is clearly $Z$, the isomorphism \eqref{E:Ziso1} follows.
It is easy to see that the isomorphism \eqref{E:Ziso1} is just the restriction to $S^G$ of
the Hopf counit $S \to \cC$, $s \mapsto s(1)$.
In particular, for any $s \in S^G$, we have 
\[
(\e\circ (1_Z \otimes \rho_l(g)))(s) = s(g^{-1}) = g.s(1) = (\rho(g) \circ \e)(s)\ ,
\]
proving \eqref{E:intertwining}.

\subsection{The case of an algebraic torus}

Now let $G \cong (\k^\times)^d$ be an algebraic torus over $\k$ and let $\Lambda = X(G)
\cong \ZZ^d$ be its lattice of rational characters. Then $\k[G] = \k\Lambda$, the
group algebra of $\Lambda$ over $\k$. As it is customary to use additive notation for the lattice
$\Lambda$, we will write the standard $\k$-basis of $\k[G]$ as $\{ \bx^\lambda \mid \lambda \in
\Lambda \}$; so $\bx^\lambda \bx^{\lambda'} = \bx^{\lambda + \lambda'}$ and
$\bx^\lambda(g) = \langle \lambda,g\rangle \in \k^\times$ for $g \in G$. Then 
$\rho_r(g)\bx^\lambda = \langle \lambda,g\rangle \bx^\lambda$ and
\[
S = \bigoplus_{\lambda \in \Lambda} \cC \otimes_\k \k \bx^\lambda \cong \cC\Lambda\ ,
\]
the group algebra of $\Lambda$ over the field $\cC$. Consider an element 
$s = \sum_\lambda s_\lambda \otimes \bx^\lambda \in S$ with $s_\lambda \in \cC$.
Then $g.s = \sum_\lambda g.s_\lambda \otimes \langle \lambda,g\rangle\bx^\lambda$
for $g \in G$.
Hence, $s \in S^G$ if and only if $g.s_\lambda = \langle -\lambda,g\rangle s_\lambda$
for all $g, \lambda$. Putting 
$\cC_\lambda = \{ c \in \cC \mid g.c = \langle \lambda,g\rangle c \text{ for all } g \in G \}$
and noting that each nonzero $\cC_\lambda$ is $1$-dimensional over the
fixed field $\cC^G$, we
have
\begin{equation} \label{E:S^G}
S^G = \bigoplus_{\lambda \in \Lambda} \cC_{-\lambda} \otimes_\k \k \bx^\lambda 
\cong \cC^G\Gamma\ ,
\end{equation}
the group algebra of the sublattice $\Gamma = \{ \lambda \in \Lambda \mid \cC_\lambda \neq 0 \}$
over $\cC^G$. We remark that $\Gamma = X(G/N)$ is the character lattice of the torus $G/N$,
where $N$ is the kernel of the action of $G$ on $Z$; so
$Z = \bigoplus_{\lambda \in \Gamma} \cC_{\lambda}$.
The isomorphisms \eqref{E:S^G} and \eqref{E:Ziso1} prove part (a) of the Theorem.
Note also the $S$ is free over $S^G$.

We claim that each $G$-stable ideal $\fa$ of $S$ is generated by its intersection
with $S^G$:
\begin{equation} \label{E:ideal}
\fa = (\fa \cap S^G)S \ .
\end{equation}
For the nontrivial inclusion $\subseteq$, let
$s = \sum_\lambda s_\lambda \otimes \bx^\lambda \in \fa$ be given.
In order to show that $s \in (\fa \cap S^G)S$, we argue 
by induction on the size of 
$\Supp(s) = \{ \lambda \in \Lambda \mid s_\lambda \neq 0 \}$, 
the length of $s$. Our claim
being clear for $s=0$, assume that $s \neq 0$. Suppose there exists an
element $0 \neq t \in \fa$ with $\Supp(t) \subsetneqq \Supp(s)$. Multiplying $t$ and $s$
with suitable units of the form $c \otimes \bx^\mu$, we may assume that
$0 \in \Supp(t)$ and $t_0=s_0=1$. Since $t$ and $s-t$ are shorter than $s$,  
they both belong to $(\fa \cap S^G)S$ and hence $s \in (\fa \cap S^G)S$ as well.
Therefore, we may assume that if $t \in \fa$ satisfies $\Supp(t) \subsetneqq \Supp(s)$
then $t=0$. Continuing to assume that $s_0=1$, this holds in particular for 
$t = s - g.s = \sum_{0 \neq \lambda} (s_\lambda - \langle \lambda,g\rangle g.s_\lambda) \otimes \bx^\lambda$
for each $g \in G$. Therefore, we must have $s \in S^G$ and \eqref{E:ideal} is proved.

Now let $\fb$ be an ideal of $S^G$ and let $\fa$ denote sum of all ideals of $S$ that contract 
to $\fb$. Since $S$ is free over $S^G$, we have $\fa \cap S^G = \fb$. Moreover, $\fa$ is clearly
$G$-stable and so \eqref{E:ideal} gives that $\fa = \fb S$. Thus, $\fb S$ is the unique largest ideal of $S$
that contracts to $\fb$.

\subsection{The prime correspondence}

We start with some reminders from \cite{mL09}. For now, let $G$ again be an arbitrary
connected affine algebraic group over $\k$ and let $\k(G) = \Fract \k[G]$ be the field of rational functions of $G$. 
The $G$-action on $S$ via $\rho \otimes \rho_r$
extends uniquely to an action of $G$ on the following localization of $S$:
\[
T = \cC \otimes_\k \k(G)\ .
\]
Let $\Spec^G(T)$ denote the collection of all $G$-stable prime ideals of
$T$. Then \cite[Theorem 9]{mL09} establishes an order isomorphism 
\begin{equation} \label{E:bij0}
c \colon \Spec_0R \stackrel{\sim}{\tto} \Spec^G(T)
\end{equation}
with the following $G$-equivariance property, for $P \in \Spec_0R$ and $g \in
G$:
\begin{equation} \label{E:equivariance}
c(g.P) = (1_\cC \otimes \rho_\ell(g))(c(P))\,.
\end{equation}

Since $T$ is the localization $S$ at the nonzero elements of $\k[G]$, 
contraction and extension yields
a $G$-equivariant order isomorphism 
$\Spec T \stackrel{\sim}{\tto} \{ \fp \in \Spec S \mid \fp \cap \k[G] = 0\}$.
Note that $\k[G]$ is a $G$-simple ring, because $G$ acts transitively on itself 
by right multiplication.
Therefore, each $\fp \in \Spec^G(S)$ satisfies $\fp \cap \k[G] = 0$, and hence 
the above bijection restricts to a bijection 
\begin{equation} \label{E:bij1}
\Spec^G(T) \stackrel{\sim}{\tto} \Spec^G(S)\ ,
\end{equation}
given by contraction and extension.

We now return to the case of an algebraic torus $G$.
Then the map $\fp \mapsto \fp \cap S^G$ injects $\Spec^G(S)$ into
$\Spec(S^G)$ by \eqref{E:ideal}. Moreover, for any $\fq \in \Spec(S^G)$, we know
that $\fp = \fq S$ is the unique largest ideal of $S$ that contracts to $\fq$, which implies that
$\fp \in \Spec^G(S)$. Thus we obtain a bijection
\begin{equation} \label{E:bij2}
\Spec^G(S) \stackrel{\sim}{\tto} \Spec(S^G)
\end{equation}
that is again given by contraction and extension.
From \eqref{E:bij0} - \eqref{E:bij2} in conjunction with the isomorphism 
\eqref{E:Ziso1} we obtain the desired order isomorphism
$\gamma \colon \Spec_0R \stackrel{\sim}{\tto} \Spec(Z)$.
The fact that $\gamma$ is $G$-equivariant is immediate from 
the $G$-equivariance property \eqref{E:equivariance}
of the bijection \eqref{E:bij0} and the intertwining formula \eqref{E:intertwining}.
This completes the proof of the Stratification Theorem.


\bibliographystyle{amsplain}
\bibliography{../bibliography}


\end{document}